\documentclass[12pt]{iopart}
\eqnobysec
\usepackage{iopams}

\usepackage{graphicx}

\newtheorem{theorem}{Theorem}

\newcommand{\diff}{{\mathrm{d}}}

\begin{document}

\title{Optimal control of two coupled spinning particles in the Euler-Lagrange picture}

\author{M. Delgado-T\'ellez$^{1}$, A. Ibort$^{2,3}$, T. Rodr\'iguez de la Pe\~na$^2$, R. Salmoni$^2$}

\address{$^{1}$Depto. de Matem\'atica Aplicada,  Univ. Polit\'ecnica de Madrid.
Avda. Juan de Herrera 6, 28040 Madrid, Spain.}
\address{$^2$Department of Mathematics, Universidad Carlos III de Madrid, Av. Universidad 30
28911 Legan\'es, Madrid, Spain.}
\address{$^3$Instituto de Ciencias Matem\'{a}ticas (CSIC - UAM - UC3M - UCM), Nicol\'{a}s Cabrera 13-15, Campus de Cantoblanco, UAM, 28049, Madrid, Spain.}

\ead{marina.delgado@upm.es, albertoi@math.uc3m.es, thrdelape@gmail.com, salmoni@math.uc3m.es}

\begin{abstract}
A family of optimal control problems for a single and two coupled spinning particles in the Euler-Lagrange formalism is discussed.   A characteristic of such problems is that the equations controlling the system are implicit and a reduction procedure to deal with them must be carried on. 

The  reduction of the implicit control equations arising in these problems will be discussed in the slightly more general setting of implicit equations defined by invariant one-forms on Lie groups.   As an instance, the first order differential equations describing the extremal solutions of an optimal control problem for a single spinning particle, obtained by using Pontryagin's Maximum Principle (PMP), will be found and shown to be completely integrable.

Then, using again PMP, solutions for the problem of two coupled spinning particles will be characterised as solutions of a system of coupled non-linear matrix differential equations.  The reduction of the implicit system will show that the reduced space for them is the product of the space of states for the independent systems, implying the absence of `entanglement' in this instance.

Finally it will be shown that, in the case of identical systems, the degree three matrix polynomial differential equations determined by the optimal feedback law,  constitute a completely integrable Hamiltonian system and some of its solutions are described explicitly.
\end{abstract}

\vspace{2pc}

\submitto{\JPA}
\maketitle

\section{Introduction}

In this paper a family of optimal control problems for single and coupled spinning particles in the Euler-Lagrange picture are discussed.  Besides, these problems could just be considered as a new family of abstract integrable optimal control problems defined on groups.

Quantum control (optimal or not) of coupled (or standalone) quantum spin systems is a relevant question in quantum control and quantum information theory and is becoming more and more relevant because of their experimental implementation achievements.  A basic requirement for quantum information processing systems is the ability to control the state of a single qubit.  Notice that for qubits based on spin, a universal single-qubit gate is realized by a rotation of the spin by any angle about an arbitrary axis \cite{Pr08}.  Moreover the control and high-fidelity readout of a nuclear spin qubit was shown in \cite{Pl13}.

Even more, individual spins, associated with vacancies in a silicon carbide lattice, have been observed and coherently manipulated \cite{Mo15}.  In the same vein electrical control of a long-lived spin qubit in a Si/SiGe quantum dot has been shown recently \cite{Ma14}.  However even if nanofabricated quantum bits permit large-scale integration, they usually suffer from short coherence times due to interactions with their solid-state environment. The outstanding challenge is to engineer the environment so that it minimally affects the qubit, but still allows qubit control \cite{Ka14}.

Let us emphasize that interacting or coupled spin systems are fundamental in quantum computation as a network of interacting and controllable spin qubits can act as a quantum computer. However, because of their magnetic and quantum-mechanical nature, the spin qubits must be controlled and measured using radically different techniques as compared to classical, transistor-based bits.  Further developments will aim at measuring and controlling the exchange interaction between pairs of spins, to demonstrate a fully functional 2-qubit quantum logic gate (see for instance the analysis of continuous feedback control in \cite{Wi12}).

Geometrical control theory has provided the mathematical background to deal with quantum spin control.  Khaneja \textit{et al} showed how to obtain efficient RF pulse trains for two-spin and three-spin NMR systems by finding sub-Riemannian geodesics on a quotient space of $SU(4)$ \cite{Kh02} and the subsequent numerical implementations of it \cite{Kh05}.

We should also mention \cite{Mo04} for a geometric control study of quantum spin systems and
\cite{Sc05} for an optimal control discussion of blocks of quantum algorithms (see also \cite[Chps. 5,6]{Al08} and the recent review of geometric optimal control for quantum systems in NMR by Bonnard \textit{et al} \cite{Bo12} and references therein).

In spite of all these developments in geometric  quantum control of spin systems, little attention has been drawn to their Euler-Lagrange picture.  The reason for this could lie in the singular nature of the Lagrangian functions describing them.  Variational descriptions for spinning particles provide a natural framework to obtain a better understanding of the corresponding quantum systems.    

There is already a long history on the Lagrangian description of the equations of motion  describing the motion of a spinning particle in the presence of an external weak homogeneous electromagnetic field \cite{Ba59} that we will not try to reproduce here (see for instance \cite{Fr96} for part of this history).   Let us just mention, because the formalism used there is close to the one that will be used in this paper, the early attempts to provide a Lagrangian picture by P. Horvathy's \cite{Ho79} and the more elaborated Lagrangian descriptions of charged particles with spin by Skagerstam and Stern \cite{Sk81} (see also a more recent discussion by Grassberger \cite{Gr01}).  

It is well-known that the Euler-Lagrange picture of a quantum spin system is not of mechanical type \cite{Ba83} and it is given by a degenerate Lagrangian function on the tangent bundle of the group $SU(2)$, i.e., the Legendre transform is not invertible.   To establish the equations of motion of the system  requires, in general, a careful analysis and it may also happen that the misleading simplicity of the quantum formulation of the problem dismisses the relevance of the analysis from such perspective.

We feel that optimal control problems of spinning systems in the Euler-Lagrange formalism deserve to be analyzed because they would definitely help in building more intuition on the behaviour of more complicated situations.    Studying such problems in the Euler-Lagrange formalism would bring together the geometrical analysis proper of optimal control problems with the geometrical picture of spinning particles.    Thus the analysis of such relevant aspects of symmetries, reduction, etc., can be done from a unified perspective and drawing significant results becomes easier as it will be shown afterwards.   

More specifically, in this paper we will analyze the optimal control problem of two coupled spinning particles in a uniform magnetic field with objective functional combining the intensity of the field plus the intensity of the coupling.    The coupled spins will be described in the Euler-Lagrange formalism as a system on the product group $SU(2)\times SU(2)$,  and it will be shown that both, the optimal control problem for a single and coupled spinning particles, define completely integrable Hamiltonian systems, offering a new insight into the structure of the corresponding quantum systems.  
 
The optimal control problems discussed in this paper will be presented in a slightly more general context, that of first order Lagrangian systems defined on Lie groups that will be discussed in Sect. \ref{sec:extremal}.  The implementation of Pontryagin's Maximum Principle (PMP) for this situation will be also considered and in the regular situation, i.e., when there exists an optimal feedback law, the reduced Hamiltonian equations satisfied by normal extremals will be derived. This will be the content of Sect. \ref{sec:reduction}.   These ideas will be applied to the discussion of two coupled spinning systems in Sect. \ref{sec:coupled}.  Finally, in Sects. \ref{sec:PMPcoupled} and \ref{sec:integrability}, the corresponding Hamiltonian equations obtained applying Pontryagin's Maximum Principle will be shown to be completely integrable and its solutions will be described using an appropriate system of coordinates.


\section{The Euler-Lagrange description of spinning particles}

We will adopt here and in what follows, the formulation used in the monograph by Balachandran \textit{et al} \cite{Ba83}  for the variational description of spinning particles.  Thus a spinning particle with spin $\mathbf{S}$ moving on a fixed external magnetic field $\mathbf{B}$ can be described by the Lagrangian function on the tangent bundle of the configuration space $\mathcal{Q} = SU(2)$, $U \in SU(2)$ (we will be discarding here the degrees of freedom corresponding to the position of the particles in space), given as:
\begin{equation}\label{lagr}
\mathcal{L}(U,\dot{U}) = i\lambda \mathrm{Tr}(\sigma_3 U^\dagger \dot{U}) + \mu \mathrm{Tr}(S B),
\end{equation}
where $\mathbf{S} = (S_1,S_2,S_3) \in \mathbb{R}^3$ is a unitary vector, $ \mathbf{S}^2=1$, and the Hermitean matrices $S$ and $B$ are defined as
$$
S = \mathbf{S}\cdot \boldsymbol{ \sigma}\, , \qquad  B= \mathbf{B} \cdot \boldsymbol{\sigma} \, ,
$$
with $\boldsymbol{\sigma} = (\sigma_1, \sigma_2, \sigma_3)$ the vector whose components are the standard Pauli matrices:
$$
\sigma_1 = \left( \begin{array}{cc} 0 & 1 \\ 1 & 0 \end{array} \right) \, , \qquad
\sigma_2 = \left( \begin{array}{cr} 0 & -i \\ i & 0 \end{array} \right) \, \qquad
\sigma_3 = \left( \begin{array}{cr} 1 & 0 \\ 0 & -1 \end{array} \right) \, .
$$

The relation between
the configuration variable $U \in SU(2)$ and the spin matrix $S$ is given by the Hopf projection map: $SU(2) \to S^2$, $U \mapsto S$:
$$
S = U^\dagger \sigma_3U  \, .
$$
In other words, if we parametrize the matrix $U$ by two complex number $z_1, z_2$ as
$$
U = U(z_1, z_2) =\left(  \begin{array}{rc}  z_1 & z_2 \\ -\bar{z}_2 & \bar{z}_1 \end{array}\right) \,
$$
satisfying $|z_1|^2 + |z_2|^2 = 1$, we will obtain:
$$
S_1 = \bar{z}_1 z_2 + \bar{z}_2 z_1 \, , \qquad S_2 = i\bar{z}_1 z_2 - i\bar{z}_2 z_1 \, , \qquad S_3 = \bar{z}_1 z_1 - \bar{z}_2 z_2 \, .
$$

The constant $\mu$ represents the magnetic moment of the system and $\lambda$ measures the spin length.

It is immediate to check that the  Euler-Lagrange equations of such system determines an implicit system of differential equations because of the Lagrangian linear dependence on the ``velocities" $\dot{U}$ of the system.  An appropriate treatment of them, based for instance on Dirac's theory of constraints (see for instance \cite{Di49}) or, in modern terms, using the Lagrangian version of the presymplectic constraints algorithm \cite{Go78}, will lead to the equations of motion of the system in Hamiltonian form:
\begin{equation}\label{evolxpS}
 \dot{S}_i= \mu \epsilon_{i j k} B_jS_k \, , \quad i = 1,2,3 \, ,
\end{equation}
where summation over repeated indices is understood (see the discussion of the inverse problem for Wong's equations in \cite{Ca95} of which the previous equations are a particular instance).

The spatial part of the system will not be considered (we are assuming an uniform magnetic field) and we will concentrate just on the spin part.   Thus the equation describing spin evolution in eq. (\ref{evolxpS}), is written in matrix notation as:
\begin{equation}\label{Smatrix}
\dot{S} = \frac{i\mu}{2} [B, S] \, ,
\end{equation}
with $B$ representing an uniform and constant magnetic field. Equation (\ref{Smatrix}), can be integrated easily to give (see \cite[Ch.1.2]{Ca14}):
$$
S(t) = e^{i\mu tB/2} S(0) e^{-i\mu tB/2} \, .
$$

Let us consider now two coupled classical spin systems.  We will follow here the prescription for composing systems in the Lagrangian formalism, that is, the configuration space of two Lagrangian systems with configuration spaces $Q_a$, $a = 1,2$, is given by $Q_! \times Q_2$.   Thus, the configuration space of the composite system will be the product group $SU(2)\times SU(2)$.   Thus both, the departing configuration space and the possible interactions among its components will be clearly shown and, after the corresponding constraint analysis, both the reduced state space of the composite system and the corresponding equations of motion will be obtained (see Section \ref{sec:discussion} for a discussion on other possibilities).    

Thus the total Lagrangian of the coupled system as a function on $T(SU(2) \times SU(2))$  will depend on pairs of unitary matrices $U_1, U_2$ and the corresponding generalized velocities $\dot{U}_1$, $\dot{U}_2$, and it will have the form:
\begin{equation}\label{total_L}
\mathcal{L} = \mathcal{L}_1 + \mathcal{L}_2 + \mathcal{L}_I \, ,
\end{equation}
where,
\begin{equation}\label{lagrangian_spin}
\mathcal{L}_\alpha (U_\alpha, \dot{U}_\alpha) = \lambda_\alpha \mathrm{Tr} (\sigma_3 U_\alpha^\dagger \dot{U_\alpha})  +\mu_\alpha \mathrm{Tr} (S_\alpha B), \qquad \alpha = 1,2 \, ,
\end{equation}
are the Lagrangians of the individual spin sytems and,
$$ \mathcal{L}_I (U_1, U_2) = \mathrm{Tr}(S_1KS_2), $$
defines the interaction between them, where the matrix $K$ determines the structure of the coupling between the spin variables.

The magnetic moments $\mu_1$, $\mu_2$ and the modules  $\lambda_1$ and $\lambda_2$ of the individual spins can be different.  Again the Euler-Lagrange equations of the system are implicit and an adaptation of Dirac's constraint algorithm must be used to determine them.  We will discuss such procedure in Section \ref{sec:reduction}.

In addition to this, a given initial configuration $(U_1,U_2)$ of the system (\ref{total_L}) can be driven by letting the parameters of the problem evolve in time: for instance, both the external magnetic field $B$ or the coupling matrix $K$ can be varied.   Thus we are led to consider the control problem determined by the Euler-Lagrange equations defined by the Lagrangian function (\ref{total_L}) with control parameters $B$ and $K$.

Moreover a natural optimal control problems for such system can be posed by  considering an objective functional $J$ depending on the variables  $U$, $K$ and $B$, that can have, for instance, the simple form (the coefficients are chosen for convenience):
\begin{equation}\label{objective}
J(K,B)=\frac{1}{2}\int_0^T ||B(t)||^2+ ||K(t)||^2  \,  dt.
\end{equation}
Thus, given an initial $U_0=(U_{1,0} , U_{2,0})$, and a target configuration $U_T=(U_{1,T},U_{2,T})$ to be reached at time $T$, we would like to know if there exist admisible curves $K(t)$, $B(t)$, and a solution $U(t)$ of the Euler-Lagrange equations
defined by (\ref{total_L}) such that $U(0)=U_0$ and $U(T)=U_T$, minimizing $J$.

Because of the implicit character of the control equation determined by the Lagrangian function $\mathcal{L}$, it is not possible to apply Pontryagin's Maximum Principle \cite{Po62} directly to this problem.   
Unfortunately a general discussion on the PMP for optimal control problems with implicit control equations is not yet available.  Some results on this direction for linear quadratic systems can be found for instance in \cite{De09} and references therein.   Other ideas involving an extension of Dirac's constraints algorithm to singular optimal problems was also discussed in \cite{delg02}, \cite{Lo00}.  An adaptation of PMP to a class of implicit differential control equations by Petit \cite{Pe98} can also be used.  However because the proposed objective functional eq. (\ref{objective}) does not depend on the state variables $(U_1,U_2)$, we can just obtain the reduced equations of motion for the Lagrangian (\ref{total_L}) and then, apply PMP to them to obtain an explicit characterization of the normal extremals of the optimal control problem. These ideas will be discussed in the following sections.
 

\section{Extremal solutions for classical spinning particles}\label{sec:extremal}

\subsection{Constraints analysis of the implicit Euler-Lagrange equations on Lie groups defined by first-order Lagrangians}\label{sec:reduction}

Because the configuration space of individual systems are groups, it is convenient to introduce a slightly more general setting that will prove to be helpful in computing the corresponding reduced systems.  Let $G$ be a Lie group that will be considered as the configuration space of a Lagrangian system. As customary $\mathfrak{g}$  will denote the Lie algebra of $G$ and $\mathfrak{g}^*$ its dual space.   Given an element $\nu \in \mathfrak{g}^*$, we denote by $\alpha_\nu$ the unique left--invariant 1--form on $G$ whose value at the identity element is $\nu$, that is $\alpha_\nu (g) = TL_{g^{-1}}^*\nu$, where $L_g$ denotes the standard left-translation by the element $g\in G$ and $TL_g$ the corresponding tangent map.    We may use left-translations to identify $TG$ with the Cartesian product $G \times \mathfrak{g}$ by means of the diffeomorphism $\Lambda \colon TG \to G \times \mathfrak{g}$, $\Lambda(g,\dot{g}) = (g, TL_{g^-1}(\dot{g}))$. Then the vertical part of the tangent bundle results canonically identified with $\mathfrak{g}$.

For instance, if $G$ is a matrix group like $SU(2)$, its elements are matrices $U = (U_{ij})$, then its Lie algebra $\mathfrak{g}$ is the linear space of matrices $A = \frac{d}{ds} U(s)\mid_{s = 0}$ where $U(s)$ is a smooth curve on $G$ such that $U(0) = I$ is the identity matrix.  Clearly in the case of the group $SU(2)$, its Lie algebra $\mathfrak{su}(2)$ is given by $2\times 2$ skew-Hermitean matrices.     The tangent space to the group $G$ at the matrix $U$ consists on matrices of the form $\dot{U} = U A$, where $A\in \mathfrak{g}$ and left (right) translation by the matrix $U^{-1}$ is given by $TL_{U^{-1}} \dot{U}= U^{-1}U A = A$.   So, the natural identification between $TG$ and $G \times \mathfrak{g}$ above is spelled as $\Lambda \colon (U, \dot{U}) \mapsto (U, A = U^{-1}\dot{U})$.

Let $\Theta^L$ be the canonical left-invariant $\mathfrak{g}$-valued Maurer--Cartan 1-form on $G$, that is    
$\Theta^L_g(\dot{g}) = TL_{g^{-1}} (\dot{g})$, $\dot{g}\in T_gG$ and,  clearly,
\begin{equation}\label{alphanu}
\alpha_\nu = \langle \nu, \Theta^L \rangle \, ,
\end{equation}
where $\langle \cdot, \cdot \rangle$ denotes the canonical pairing between $\mathfrak{g}$ and its dual space $\mathfrak{g}^*$.   Because of the Maurer-Cartan equation
\begin{equation}\label{MC}
\diff \Theta^L + \frac{1}{2} \Theta^L \wedge \Theta^L = 0\, ,
\end{equation}
it is also evident that:
\begin{equation}\label{dalpha}
\diff \alpha_\nu = - \frac{1}{2} \langle \nu, \Theta^L \wedge \Theta^L \rangle \, .
\end{equation}

Again, in the particular instance of matrix groups, the previous definitions become particularly simple.  Thus for instance the canonical left-invariant Cartan 1-form becomes the matrix valued 1-form whose $(i,j)$ component is given by  $\Theta^L_{ij} = U^{-1}_{ik} dU_{kj}$.  However in general, it is more convenient to choose a linear basis $\zeta_a$, $a = 1, \ldots, \dim \mathfrak{g}$ and its dual basis $\theta^a$, $\langle \theta^a, \zeta_b \rangle = \delta_a^b$.  Then consider the 1-family of left invariant vector fields on $G$ defined by the elements $\zeta_a$ (and the corresponding left-invariant 1-forms too) and use them to write explicit formulas for $\Theta^L$ and other geometrical objects.  Denoting with the same symbol the element $\zeta_a$ and the corresponding left-invariant vector field, and doing the same for $\theta^a$, we get easily that $\Theta^L = \zeta_a \otimes \theta^a$.    Thus if $\alpha_\nu$ denotes the left-invariant 1-form whose value at the identity is $\nu$, then $\alpha_\nu = \nu_a \theta^a = \langle \nu, \Theta^L \rangle$ as in Eq. (\ref{alphanu}).  With these notations Maurer-Cartan equations (\ref{MC}) become simply:
$$
d\theta^a = - C_{bc}^a \, \theta^b \wedge \theta^c \, ,
$$
where $C_{bc}^a$ denote the structure constants of the Lie algebra $\mathfrak{g}$ in the basis $\zeta_a$, i.e., $[\zeta_b \zeta_c] = C_{bc}^a \zeta_a$, 
and Eq. (\ref{dalpha}) reads:  
\begin{equation}\label{dalphacomp}
d\alpha_\nu = - \nu_a C_{bc}^a  \theta^b \wedge \theta^c \, .
\end{equation} 

Finally, notice that for semisimple compact groups (like $SU(2)$) the Killing-Cartan form is negative non-degenerate and it allows to identify $\mathfrak{g}$ and its dual $\mathfrak{g}^*$, hence the bracket $\langle\cdot, \cdot \rangle$ denotes either the canonical pairing between $\mathfrak{g}$ and its dual $\mathfrak{g}^*$ or the Killing-Cartan form.  In the particular instance of $SU(n)$ groups, the Killing-Cartan form is given as  $\langle A, B\rangle = - \frac{1}{2} \mathrm{Tr} (A^\dagger B)$ and we may identify $\mathfrak{su}(n)$ naturally with its dual space $\mathfrak{su}(n)^*$.  Again in the particular instance of $SU(2)$, once we choose and orthonormal basis for $\mathfrak{su}(2)$, we may identify it with $\mathbb{R}^3$ and its Euclidean metric.   This correspondence is truly what lies at the bottom of the identification between vectors $\mathbf{S}$ in $\mathbb{R}^3$ and Hermitean matrices $S$ used so far (properly speaking the identification is between vectors $\mathbf{S}$ and skew-Hermitean matrices $\hat{S} = - \frac{i}{2}S$ as it will be explained below, Sect. \ref{sec:optimal}).

Now, either using Eq. (\ref{dalpha}) or its components expression Eq. (\ref{dalphacomp}), it is simple to check that given a left-invariant 1-form $\alpha_\nu$, the characteristic distribution $K= \ker \omega_\nu$ of the presymplectic form $\omega_\nu = \diff\alpha_\nu$ on $G$, is given by $K = \mathfrak{g}_\nu$, where $\mathfrak{g}_\nu$ denotes the isotropy algebra of $\nu$ with respect to the coadjoint action, i.e., the Lie algebra of the isotropy group $G_\nu = \{ g \in G \mid \mathrm{Ad}_g^* \nu = \nu \}$, more explicitly 
\begin{equation}\label{gnu}
K = \mathfrak{g}_\nu = \{Ê\xi \in \mathfrak{g} \mid \langle \nu , [\xi, \zeta] \rangle = 0 \, , \forall \zeta \in \mathfrak{g} \} \, .
\end{equation}

The characteristic distribution $\ker \omega_\nu$ is integrable because $\omega_\nu$ is closed, and the connected components of the leaves of the foliation $\mathcal{K}$ defined by it, are orbits of the left action of the isotropy group $G_\nu$ on $G$.

Consider now the Lagrangian system on $G$ with Lagrangian function $\mathcal{L} \colon TG \to \mathbb{R}$ given by:
\begin{equation}\label{L}
\mathcal{L} (g,\dot{g}) = \langle \alpha_\nu(g) , \dot{g} \rangle - V(g)\, , \qquad (g,\dot{g}) \in TG ,
\end{equation}
with $V\colon G \to \mathbb{R}$ a $G_\nu$--invariant function on $G$.  The Poincar\'e--Cartan 1--form $\theta_{\mathcal{L}} = \partial \mathcal{L} / \partial{\dot{g}} dg$ of the system is easily obtained to be:
\begin{equation}\label{thetaL}
\theta_\mathcal{L} = \tau^* \alpha_\nu \, ,
\end{equation}
where $\tau \colon TG \to G$ denotes the canonical projection of the tangent bundle of $G$.    Then, the Cartan 2--form $\omega_{\mathcal{L} } = - d\theta_{\mathcal{L}}$ is just:
\begin{equation}\label{omegaL}
\omega_{\mathcal{L}} = - \tau^* d\alpha_\nu = \frac{1}{2} \langle  \nu, \tau^*( \Theta^L \wedge \Theta^L )\rangle \, ,
\end{equation}
and the kinetic energy of the system is simply given by
\begin{equation}\label{EL}
E_{\mathcal{L}} = \dot{g} \frac{\partial \mathcal{L}}{\dot{g}} - \mathcal{L} = V \, .
\end{equation}

Hence the Euler-Lagrange vector field $\Gamma$ for the Lagrangian function eq. (\ref{L}) are given by the implicit system of equations on $TG$,
$$
i_\Gamma \omega_\mathcal{L} = \diff E_\mathcal{L} \, ,
$$
which, because of eqs. (\ref{thetaL}), (\ref{omegaL}) and (\ref{EL}), are equivalent to:
\begin{equation}\label{ELV}
i_\Gamma  (\tau^* d\alpha_\nu) = -\tau^* \diff V \, .
\end{equation}

We apply now the constraint algorithm to eqs. (\ref{ELV}) \cite{Go78}.   We have to characterize first the characteristic distribution of $\omega_{\mathcal{L}}$.  Because of the previous discussion is easily seen that
$$
\ker \omega_{\mathcal{L}} = \mathfrak{g}_\nu \oplus \mathfrak{g} \, ,
$$
where we have used the identification $TG \cong G \times \mathfrak{g}$ discussed above.

Notice that the Hamiltonian of this system, Eq. (\ref{ELV}), is $\tau^*V$ which is invariant with respect to the vector fields on $\ker \omega_{\mathcal{L}}$
because by definition it is $G_\nu$-invariant.  Consequently the constraints algorithm stops at the first step and the final constraints submanifold is simply the total space $TG$.
However the dynamical equations (\ref{ELV}) have a large kernel and the reduced state space of the system is given by $TG/\ker \omega_{\mathcal{L}} $.

This quotient can be computed in two steps.  First we will quotient the system with respect to the vertical distribution $\mathfrak{g}$ and secondly with respect to the characteristic distribution $\mathfrak{g}_\nu$ (which is tangent to the configuration space $G$ considered as the zero section of the tangent bundle $TG$).   In fact because of the identification $TG \cong G \times \mathfrak{g}$, it is obvious that $TG/\mathfrak{g} \cong G$ and, after the first step,  the reduced system is defined on $G$ and takes the simple form:
\begin{equation}\label{first_red}
i_\Gamma \diff\alpha_\nu = \diff V .
 \end{equation}

This system still has a kernel $K\ker d\alpha_\nu$, i.e., it defines a presymplectic Hamiltonian system on $G$.   Now  because of the discussion before and after Eq. (\ref{gnu}), we get that $K = \mathfrak{g}_\nu$, and the true reduced space happens to be $G/G_\nu \cong \mathcal{O}_\nu$ where $\mathcal{O}_\nu$ denotes the coadjoint orbit of $G$ passing through $\nu$.  In fact, a simple computation shows that the projection of the presymplectic form $\diff \alpha_\nu$ to $G/G_\nu$ gives the canonical Kostant-Kirillov-Souriau symplectic structure on $\mathcal{O}_\nu$.

We must notice that if $V$ were not $G_\nu$--invariant the constraint algorithm would have to be pursued until obtaining the final constraint submanifold of the problem.

\subsection{An optimal control problem for a single spinning system}\label{sec:optimal}

We will prepare the ground for the study of two coupled spinning systems by considering first the case of
the optimal control of a single classical spinning particle consisting on reaching a prescribed state $S_1$ starting from a given one $S_0$ in a fixed time $T$.  In such case the configuration space of the system will be the group  $SU(2)$ and its (un)reduced state space will be $TSU(2)$.  The equations of the system will be given implicitly by the Lagrangian function (\ref{lagr}) and,  comparing with the general form of a Lagrangian defined by a left-invariant 1-form on a Lie group $G$, eq. (\ref{L}), we have that $\nu = -\frac{i}{2}\lambda \sigma_3$ and $V(U) = -\mu \mathrm{Tr}(SB)$.

We will consider now the following simple objective functional (the coefficient is chosen for convenience):
\begin{equation}\label{J0}
J_0 ( B) = -\frac{1}{8} \int_0^T  \mathrm{Tr} (B(t)^2) \,  dt \,,
\end{equation}
as there are no coupling term in this situation.

The extremal solutions of the optimal control problem defined above, are obtained by solving the optimal control problem defined by the objective functional eq. (\ref{J0}) restricted to the the reduced state space of the system.   In this case because of the analysis in the previous section, Sect. \ref{sec:reduction}, such reduced space is $S^2$, which is the quotient of $SU(2)$ with respect to the characteristic distribution of the presymplectic form $\diff \alpha_\nu$ which is just the diagonal subgroup $U(1)$.  Hence $SU(2) / U(1) \cong S^2$ and the canonical projection $\varrho \colon SU(2) \to S^2$ is just the Hopf map.

Notice that the equation of motion are given by Eq. (\ref{first_red}) with the potential $V$ being clearly $U(1)$-invariant.   Recall that the canonical embedding of the coadjoint orbit $\mathcal{O}_\nu \subset \mathfrak{g}^*$ allows to write the reduced Hamilton's equations as the restriction to $\mathcal{O}_\nu$ of the Poisson equations defined on $\mathfrak{g}^*$ by the linear form defined by $V$ with respect to the canonical Lie-Poisson structure on it.   In our case, $\mathcal{O}_\nu \cong S^2 \subset \mathfrak{su}(2)^* \cong \mathbb{R}^3$, and we get immediately Eq. (\ref{evolxpS}) or, in matrix form, Eq. (\ref{Smatrix}).   However we will rather write the equation on motion in terms of natural variables on the Lie algebra $\mathfrak{su}(2)$, i.e., skew-hermitian matrices $\hat{S} = -\frac{i}{2} S$, $\hat{B} =  -\frac{i}{2} B$.  Notice that the matrices $\hat{S}$, $\hat{B}$ satisfy $[\hat{S}, \hat{B}] = \widehat{\mathbf{S} \wedge \mathbf{P}}$, with $\mathbf{S}$, $\mathbf{B}$ the vectors in $\mathbb{R}^3$ associated to $\hat{S}$ and $\hat{B}$ respectively.   Then, Eq. (\ref{Smatrix}) becomes:
$$
\frac{d}{dt} \hat{S} = \mu [\hat{S}, \hat{B}] \, .
$$
In what follows we will use the matrix notation $\hat{S}$ or the vector notation $\mathbf{S}$ depending on the context.

To apply PMP we construct first Pontryagin's bundle, that is the space $M = T^*S^2 \times \mathbb{R}^3$, consisting of points $(\mathbf{S},\mathbf{P}; \mathbf{B})$, with $\mathbf{S}\cdot \mathbf{S} = 1$, and $\mathbf{P}\cdot \mathbf{S} = 0$.   Using matrix notation we will have $|| \hat{S} || = 1$, $\langle \hat{P} , \hat{S} \rangle = 0$.  The variables $\mathbf{S}$ represent the states of the system, $\mathbf{P}$ are called the co-estate variables and $\mathbf{B}$ are the controls of the problem.    Pontryagin's Hamiltonian is given by:
$$
H_P (\hat{S}, \hat{P}; \hat{B}) = \langle \hat{P}, \mu [\hat{S}, \hat{B}] \rangle - \frac{1}{2} \langle \hat{B}, \hat{B} \rangle \, ,
$$
that leads to the adjoint equations for the co-estate variables $\hat{P}$:
$$
\frac{d}{dt} \hat{P} = - \frac{\partial H_P}{\partial \hat{S}} = \mu [\hat{P}, \hat{B}] \, ,
$$
and the optimal feedback law:
\begin{equation}\label{BSP}
0 = \frac{\partial H_P}{\partial \hat{B}} = \mu [\hat{P}, \hat{S}] - \hat{B} \, .
\end{equation}
that provides the magnetic field $\mathbf{B}$ that should be applied at each time to the system.

Substituting the optimal feedback relation (\ref{BSP}) in Pontryagin's Hamiltonian $H_P$, we get:
\begin{equation}\label{HSxB}
H_P (\hat{S}, \hat{P}) = \frac{1}{2} \langle \hat{B}, \hat{B} \rangle = \frac{\mu^2}{2} || [\hat{S}, \hat{P}] ||^2 \, .
\end{equation}

The Hamiltonian system on $T^*S^2$ defined by the Hamiltonian function Eq. (\ref{HSxB}) is completely integrable.
The equations of motion are given by:
\begin{equation}\label{SPmu}
\frac{d}{dt} \hat{S} = \mu^2 [\hat{S}, [\hat{P} , \hat{S}]] \, , \qquad \frac{d}{dt} \hat{P} = \mu^2 [\hat{P}, [\hat{P} , \hat{S}]] \, ,
\end{equation}
and it is easy to check that $[\hat{S}, \hat{P}]$ is a constant of the motion together with the Hamiltonian itself.
Actually we may consider the Hamiltonian system on $T^*\mathbb{R}^3$ defined by the Hamiltonian (no restrictions on the modulus of $\mathbf{S}$):
$$
H(\mathbf{S}, \mathbf{P}) = \frac{\mu^2}{2} || \mathbf{S} \wedge \mathbf{P} ||^2 \, .
$$
This system can be easily integrated by observing that $\mathbf{S}\cdot \mathbf{S}$, $\mathbf{S}\cdot \mathbf{P}$,  $\mathbf{P}\cdot \mathbf{P}$ and $\mathbf{S}\wedge \mathbf{P}$  are constants of the motion. Thus the motion takes place, given initial values $\mathbf{S}_0$ and  $\mathbf{P}_0$ for the  momenta, in the bundle of spheres of radius $|| \mathbf{P}_0 ||$ over $S^2$.   In particular, selecting the level set corresponding to $\mathbf{S}\cdot \mathbf{S} = 1$, $\mathbf{S}\cdot \mathbf{P} = 0$, the restricted system becomes (\ref{HSxB}).

The solutions of such system can be described explicitly as follows.   The magnetic field is given by Eq. (\ref{BSP}), or in vector notation $\mathbf{B} = \mu \mathbf{P} \wedge \mathbf{S}$, but  because $\mathbf{S} \wedge \mathbf{P}$ is constant of the motion then $\mathbf{B}$ is constant in time and perpendicular to both $\mathbf{S}$ and $\mathbf{P}$.  Then the motion of $\mathbf{S}$ is a rotation around $\mathbf{B}$ with angular velocity $|| \mathbf{B} ||/ \mu$. Notice that the same happens for $\mathbf{P}$ because the evolution of $\mathbf{S}$ determines that of $\mathbf{P}$.

Thus, given an initial state $\mathbf{S}_0 \in S^2$, any state $\mathbf{S}_1 \in S^2$ can be reached in time $T$ following an optimal trajectory.   For that we must pick up a tangent vector $\mathbf{P}_0$ which is coplanar with $\mathbf{S}_0, \mathbf{S}_1$, orthogonal to $\mathbf{S}_0$ and of norm $\theta_0 / \lambda T$ where $\theta_0$ denotes the angle determined by $\mathbf{S}_0$ and $\mathbf{S}_1$.   Then the magnetic field $\mathbf{B}$ is orthogonal to the plane defined by $\mathbf{S}_0$ and $\mathbf{S}_1$, and the spin variable will rotate around $\mathbf{B}$ with angular velocity $ \lambda || \mathbf{P}_0 ||$.   Notice that the vector $\mathbf{P}_0$ is uniquely defined unless $\mathbf{S}_1$ is antipodal to $\mathbf{S}_0$ in which case any tangent vector $\mathbf{P}_0$ of the appropriate length will suffice.  Hence we conclude that the system is state controllable and the trajectories joining two states can be chosen to be optimal.


\section{Optimal control of two coupled spinning particles}


\subsection{Constraints analysis of the implicit Euler-Lagrange equations for coupled spinning particles}\label{sec:coupled}

We can apply now the previous results to the problem at hand.   First we want to obtain the reduced state space of the system defined by the Lagrangian function defining two coupled spinning particles.  Such Lagrangian function is defined on the tangent bundle of $SU(2)\times SU(2)$ and has the form given in Eq. (\ref{total_L}).   The configuration variables will be pairs $(U_1,U_2)$ of $2\times 2$ special unitary matrices (in what follows we will use a subindex $\alpha = 1,2$ to label them.)

The Lagrangians $\mathcal{L}_\alpha$ of each individual system have the form given in Eq. (\ref{lagrangian_spin}), that, in the intrinsic form described in Sect. \ref{sec:reduction}, have the form of first order Lagrangian on groups,  Eq. (\ref{L}), where the chosen left-invariant 1-forms $\nu_\alpha$, $\alpha = 1,2$, in $\mathfrak{su}(2)^*$ are given by
$$
\nu_\alpha = -\frac{i}{2} \lambda_\alpha \sigma_3 \, , \qquad \alpha = 1,2 \, ,
$$
where we have used as usual the canonical Killing-Cartan form on $\mathfrak{su}(2)$ to identify it with its dual space.
Hence, the Lagrangian function corresponding to the composite system has the following form:
\begin{eqnarray}\label{LUalpha}
\mathcal{L} (U_1,U_2,\dot{U}_1, \dot{U}_2) &=&  \langle \alpha_{\nu_1}(U_1), \dot{U}_1\rangle +\langle \alpha_{\nu_1}(U_1), \dot{U}_1\rangle
+ \nonumber \\ &+& V_1(U_1) + V_2(U_2) + V_I (U_1,U_2) \, ,
\end{eqnarray}
where the potential functions $V_1$, $V_2$, $V_I$, are given respectively by the expressions,
$$
V_\alpha (U_\alpha ) = \mu_\alpha \mathrm{Tr} ( U_\alpha^\dagger \sigma_3 U_\alpha B ) \, , \qquad V_I (U_1, U_2) =  \mathrm{Tr} (U_1^\dagger \sigma_3 U_1 K  U_2^\dagger \sigma_3 U_2) \, .
$$
Hence repeating the computations leading to the Cartan 2--form of a single spin, Eq. (\ref{omegaL}), we get now:
$$
\omega_{\mathcal{L}} = \tau_1^* d\alpha_{\nu_1} + \tau_2^* d\alpha_{\nu_2} \, ,
$$
with $\tau_\alpha$, $\alpha = 1,2$, denoting the canonical projections $\tau_\alpha (U_1,U_2,\dot{U}_1, \dot{U}_2) = U_\alpha$.
The characteristic distribution of $\omega_{\mathcal{L}}$ will have the form:
$$
\ker \omega_{\mathcal{L}} = \ker \tau_1^* d\alpha_{\nu_1} \oplus
\ker \tau_2^* d\alpha_{\nu_2} \cong (\mathfrak{u}(1) \oplus \mathfrak{su}(2)  ) \oplus (\mathfrak{u}(1) \oplus \mathfrak{su}(2)) \, ,
$$
where $\mathfrak{u}(1)$ represents the Lie algebra of the isotropy group $SU(2)_{i\sigma_3}$ of the element $i\sigma_3$ in the dual of the Lie algebra $\mathfrak{su}(2)$, i.e., the one--parameter subgroup $\{ U_3 (s) = e^{is\sigma_3} \}$.

After performing the first reduction step as in the case of eq. (\ref{first_red}), we will obtain the presymplectic system on the product group $SU(2) \times SU(2)$ defined by the closed 2--form $\omega = d\alpha_{\nu_1} + d\alpha_{\nu_2}$ and Hamiltonian $V (U_1, U_2) = V_1(U_1) + V_2(U_2) + V_I (U_1, U_2)$.   All terms in the Hamiltonian are obviously invariant with respect to the isotropy subgroup $U(1) \times U(1)$ acting on the left on $SU(2)\times SU(2)$ (each one of the components $U(1)$ has the form above $U_3(s)$).  

Notice again that the subgroup $U(1) \times U(1)$ spans the characteristic distribution of $\omega$. Thus the system projects to the quotient $SU(2) \times SU(2) / (U(1) \times U(1))$ which is trivially diffeomorphic to the product of two spheres $S^2 \times S^2$.   The projection of the presymplectic form $\omega$ is given explicitly as $\frac{1}{4\pi} (\omega_{\nu_1} + \omega_{\nu_2})$, where each factor $\omega_{\nu_\alpha}$ denotes the canonical area 2--form on the sphere of radius $\lambda_\alpha$.

Using again natural spin variables
$\mathbf{S}_\alpha$, $\alpha = 1,2$, the canonical commutation relations defined by the induced symplectic structure above on the quotient space $S^2 \times S^2$, among the components $S_{\alpha i}$, $\alpha = 1,2$, $i = 1,2,3$, of the spin variables take the simple form:
\begin{equation}\label{comm_relations}
\{ S_{\alpha i}, S_{\beta j} \} = \frac{1}{\lambda_\alpha^2} \delta_{\alpha \beta}  \epsilon_{ijk} S_{\alpha j} S_{\beta k} \, , \qquad \alpha, \beta = 1,2\, , \quad i,j = 1,2,3 \, .
\end{equation}

Thus we conclude  that the final reduced space of two coupled spin systems with Lagrangian (\ref{LUalpha}) consists of the Cartesian product of the state spaces of the individual systems.   In this sense we have shown there is no ``entanglement'' in the system determined by Lagrangian (\ref{LUalpha}) as the states of the composite system are pairs of individual states even if the dynamics induced on such space is not separable, i.e., it is not a direct sum of individual dynamics because of the term $V_I$ in the Lagrangian function.    

It is important to observe that other choices for the interaction potential $V_I$ could have been considered.   For instance, it is possible to consider interaction potentials of the form $V_I (U_1, U_2) = \mathrm{Tr}(K U_1^\dagger U_2)$ which is clearly invariant only under the diagonal subgroup $U(1)$ in $SU(2) \times SU(2)$.  In such a case, the set of states such that $V_I$ would be invariant with respect the characteristic distribution of $\omega$ will reduce only to the pairs $U_1 = U_2$.   Thus, the constraints algorithm will impose that, restricted to such subspace, the interaction term will be constant and the coupled spin systems will be be trivial.

Finally we notice that if the matrix $K$ introduced in the interaction term $V_I$ is a multiple of the identity, then the system exhibits an additional symmetry corresponding to the right action of the $U(1)$ subgroup of $SU(2)$ that leaves invariant the magnetic field $B$.    The main consequence of such situation is that the system becomes integrable as it will be discussed in the next section. Thus we will assume in what follows that $K = \kappa \mathbb{I}_2$.  Then, the equations of motion on $S^2\times S^2$ given by projecting eqs. (\ref{first_red}), can be written in matrix notation as:
\begin{equation}\label{reduced}
\dot{S}_1 = \frac{i\mu_1}{2} [B, S_1] + \frac{i\kappa}{2} [S_2, S_1], \quad \dot{S}_2 =  \frac{i\mu_2}{2}  [B, S_2] +  \frac{i\kappa}{2}  [S_1,S_2].
\end{equation}
which obtained easily by computing $\diff V = \diff V_1 + \diff V_2 + \diff V_I$ and using the commutation relations (\ref{comm_relations}).


\subsection{Optimal trajectories and PMP for coupled spinning particles}\label{sec:PMPcoupled}

We are now ready to characterise the extremal trajectories of the optimal control problem  posed by the objective functional (\ref{objective}) on the system of coupled spinning particles described by the Lagrangian function (\ref{total_L}) (or, equivalently (\ref{LUalpha})).    Because of the regularity of the constraints analysis we apply Pontryagin's Maximum Principle, as we did in the case of a single spin, to the system  (\ref{reduced}) instead.

The quantities $B$ and $\kappa$ will be considered as the control variables of the system.  
 If no restrictions on the values of such variables are introduced, then PMP assures that normal extremals for the objective functional $J$, Eq. (\ref{objective}), will be given by integral curves of the Hamiltonian equations defined by Pontryagin's Hamitonian function $H_P$ defined on $T^*S^2 \times T^*S^2$  (depending on $B$ and $\kappa$ too).  
 
As in the case of a single spin system it is convenient to use natural Lie algebra variables to describe it, that is, we will consider in what follows skew-Hermitean matrices $\hat{S}_\alpha$ and $\hat{P}_\alpha$, $\alpha = 1,2$, as in Sect. \ref{sec:optimal} to denote state and co-estate variables.  Now the state Eqs. (\ref{reduced}) become:
\begin{equation}\label{reducedhat}
\frac{d}{dt}\hat{S}_1 = \mu_1 [\hat{S}_1, \hat{B}] + \kappa [\hat{S}_1, \hat{S}_2]\, , \quad \frac{d}{dt}\hat{S}_2 =  \mu_2  [\hat{S}_2, \hat{B}] +  \kappa  [\hat{S}_2,\hat{S}_1] \, .
\end{equation}

With the notations above Pontryagin's Hamitonian has the form:
\begin{eqnarray}\label{Pcoupled}
H_P (\hat{S}_1, \hat{S}_2, \hat{P}_1, \hat{P}_2; \hat{B}, \kappa ) &=& \langle \hat{P}_1, \mu_1 [\hat{B}, \hat{S}_1] + \kappa [\hat{S}_1, \hat{S}_2] \rangle + \\ &+& \langle \hat{P}_2, \mu_2 [\hat{B}, \hat{S}_2] +  \kappa  [\hat{S}_2,\hat{S}_1] \rangle - \frac{1}{2} \langle \hat{B},\hat{B} \rangle - \frac{1}{2} \kappa^2 \, , \nonumber  
\end{eqnarray}
where the co-estate variables $\hat{P}_\alpha$ denote canonical momenta in $T^*S^2$ and are such that $\langle \hat{P}_\alpha, \hat{S}_\alpha \rangle = 0$.  Notice that if $\mathbf{P}_\alpha$ denotes as usual the vector in $\mathbb{R}^3$ associated to the Hermitean matrix $P_\alpha$, then the last conditions amounts to $\mathbf{P}_\alpha \cdot \mathbf{S}_\alpha = 0$, $\alpha =  1,2$.

The regular character of the objective functional implies the existence of an optimal feedback law given by $\partial H_P /\partial \kappa = 0$, that is:
\begin{equation}\label{optimalkhat}
 \kappa =  \langle \hat{P}_1 - \hat{P}_2, [\hat{S}_1,\hat{S}_2]\rangle \, ,
 \end{equation}
and, $\partial H_P /\partial B = 0$, which is equivalent to:
\begin{equation}\label{optimalBhat}
\hat{B} = \mu_1 [\hat{P}_1,\hat{S}_1] + \mu_2  [\hat{P}_2,\hat{S}_2] \, .
\end{equation}
The adjoint equations are given by
$$
\frac{d}{dt}\hat{P}_\alpha = - \frac{\partial H_P}{\partial \hat{S}_\alpha} \, , \qquad \alpha = 1,2 \, ,
$$
that is:
\begin{eqnarray}
\frac{d}{dt}\hat{P}_1 &=& \mu_1  [\hat{P}_1,\hat{B}] + \kappa [\hat{P}_1 - \hat{P}_2 , \hat{S}_2], \nonumber \\ \frac{d}{dt}\hat{P}_2 &=& \mu_2 [\hat{P}_2,\hat{B}] + \kappa [\hat{P}_2 - \hat{P}_1,\hat{S}_1] \, . \label{Pdothat}
\end{eqnarray}


\subsection{Controllability and integrability of an optimal control problem for coupled spining particles}\label{sec:integrability}

Notice that the system of equations (\ref{reducedhat})-(\ref{Pdothat}) obtained by substituting the optimal feed-back relations Eqs. (\ref{optimalkhat})-(\ref{optimalBhat}) in them, constitute a rather complicated nonlinear coupled system of third-order matrix polynomial equations in the variables $\hat{S}_\alpha$, $\hat{P}_\alpha$ subjected to the constraints $|| \hat{S}_\alpha || = \lambda_\alpha$, and $\langle \hat{P}_\alpha, \hat{S}_\alpha \rangle = 0$, $\alpha = 1,2$.

However if we consider now the case of identical spin systems, i.e., $\mu_1 = \mu_2$ and we substitute the feedback laws Eqs. (\ref{optimalkhat})-(\ref{optimalBhat}) into the expression for Pontryagin's Hamiltonian Eq. (\ref{Pcoupled}) we obtain the simple formula:
$$
H_P (\hat{S}_1, \hat{S}_2, \hat{P}_1, \hat{P}_2; \hat{B}, \kappa ) = \frac{1}{2}\langle \hat{B}, \hat{B} \rangle  + \frac{1}{2}\langle \kappa, \kappa \rangle  \, ,
$$
or, in terms of the original variables $\hat{S}_\alpha$, $\hat{P}_\alpha$ in the phase space $T^*S^2 \times T^*S^2$, we get:
\begin{eqnarray}\label{Hreducedcoupled}
H_P (\hat{S}_1, \hat{S}_2, \hat{P}_1, \hat{P}_2; \hat{B}, \kappa ) = && \nonumber \\ = \frac{\mu^2}{2} || [\hat{P}_1,\hat{S}_1] + [\hat{P}_2,\hat{S}_2] ||^2 && + \frac{1}{2}|\langle \hat{P}_1 - \hat{P}_2, [\hat{S}_1,\hat{S}_2] \rangle |^2 \, .
\end{eqnarray}

\begin{theorem}  The Hamiltonian system on $T^*(S^2 \times S^2)$ described by the Hamiltonian function above, eq. (\ref{Hreducedcoupled}) is completely integrable.
\end{theorem}

The proof is easily obtained by realising that the quantity $ \mu [\hat{P}_1,\hat{S}_1] + \mu  [\hat{P}_2,\hat{S}_2]$ is a constant of the motion.  Actually, a long, but easy, computation shows that:
$$
\{ H_P, \mu [\hat{P}_1,\hat{S}_1] + \mu  [\hat{P}_2,\hat{S}_2] \} = 0 \, ,
$$
and
$$
\{ H_P , \kappa \} = 0 \, ,
$$
where $\{ \cdot, \cdot \} $ denotes the canonical Poisson bracket on the cotangent bundle $T^*(S^2\times S^2)$.

For instance, using now the vector notation, we may compute $\dot{\bf{B}}$ (up to a proportionality factor) as follows: 
\begin{eqnarray*}
\dot{\bf{B}}& \propto &\dot{\bf{S_1}}\wedge\bf{P_1}+\dot{\bf{S_2}}\wedge \bf{P_2}+ {\bf{S_1}}\wedge\dot{\bf{P_1}}+{\bf{S_2}}\wedge\dot{\bf{P_2}}=\\
&=&\mu(\bf{S_1}\wedge\bf{B})\wedge\bf{P_1}+\kappa(\bf{S_1}\wedge\bf{S_2})\wedge\bf{P_1}+\\
&+&\mu(\bf{S_2}\wedge\bf{B})\wedge\bf{P_2}+\kappa(\bf{S_2}\wedge\bf{S_1})\wedge\bf{P_2}+\\
&+&\mu\bf{S_1}\wedge(\bf{P_1}\wedge\bf{B})+\kappa\bf{S_1}\wedge[(\bf{P_1}-\bf{P_2})\wedge\bf{S_2}]+ \\ &+&\mu\bf{S_2}\wedge(\bf{P_2}\wedge\bf{B})
+\kappa\bf{S_2}\wedge[(\bf{P_2}-\bf{P_1})\wedge\bf{S_1}].
\end{eqnarray*}
The terms in the right hand side of the previous equation can be split in three parts:
\begin{eqnarray*}
(I)&:=&\kappa(\bf{S_1}\wedge\bf{S_2})\wedge\bf{P_1}+\kappa(\bf{S_2}\wedge\bf{S_1})\wedge\bf{P_2},\\
(II)&:=&\kappa\bf{S_1}\wedge[(\bf{P_1}-\bf{P_2})\wedge\bf{S_2}]+\kappa\bf{S_2}\wedge[(\bf{P_2}-\bf{P_1})\wedge\bf{S_1}],\\
(III)&:=&\mu(\bf{S_1}\wedge\bf{B})\wedge\bf{P_1}+\mu(\bf{S_2}\wedge\bf{B})\wedge\bf{P_2}+\mu\bf{S_1}\wedge(\bf{P_1}\wedge\bf{B})+\mu\bf{S_2}\wedge(\bf{P_2}\wedge\bf{B}).
\end{eqnarray*}

Using the vector and scalar products properties together with the constraints satisfied by the variables $\mathbf{S}_\alpha$, $\mathbf{P}_\alpha$, $\alpha = 1,2$, we get that $(II)=-(I)$ and that $(III)=0$, thus $\dot{\bf{B}}=0$.

\medskip

It is also easy to check that $\{ \hat{B} , \kappa \}= 0$.  Notice that the Hamiltonian vector field $X_\kappa$ associated to the function $\kappa$ is given by:
\begin{eqnarray}
X_\kappa &=& \frac{\partial \kappa}{\partial \hat{P}_1} \cdot \frac{\partial }{\partial \hat{S}_1} + \frac{\partial \kappa}{\partial \hat{P}_2} \cdot \frac{\partial }{\partial \hat{S}_2} - \frac{\partial \kappa}{\partial \hat{S}_1} \cdot \frac{\partial }{\partial \hat{P}_1}  - \frac{\partial \kappa}{\partial \hat{S}_2} \cdot \frac{\partial }{\partial \hat{P}_2} =\nonumber  \\
&=& ([\hat{S}_1, \hat{S}_2  ]) \cdot \frac{\partial }{\partial \hat{S}_1} - ([\hat{S}_1, \hat{S}_2  ]) \cdot \frac{\partial }{\partial \hat{S}_2} -  \nonumber \\ &&- ([\hat{S}_2,\hat{P}_1-\hat{P}_2]) \cdot \frac{\partial }{\partial \hat{P}_1}  + ([\hat{S}_1,\hat{P}_1-\hat{P}_2]) \cdot \frac{\partial }{\partial \hat{P}_2} \, .
\end{eqnarray}
Then, we compute $\{ \hat{B}, \kappa \} = X_\kappa (\hat{B}) = 0$.  

Because the Hamiltonian system has dimension 4, the constants of the motion given by the Hamiltonian itself and $\kappa$ suffice to integrate the system.  However it is possible to find a set of independent variables exhibiting a simple dependence on the spin and the corresponding co-estate variables, that will provide an explicit integration of the system.  Consider the variables:
$$
\hat{S}_\pm = \hat{S}_1 \pm \hat{S}_2 \, ; \qquad \hat{P}_\pm = \hat{P}_1 \pm \hat{P}_2 \, .
$$
We find
\begin{equation}\label{splusdot}
\frac{d}{dt} \hat{S}_+ = \mu [\hat{S}_+, \hat{B}] \, ,\qquad \frac{d}{dt} \hat{P}_- = [\hat{P}_-, \mu\hat{B} - \kappa \hat{S}_+] \, ,
\end{equation}
that together with $\hat{B}$ and $\kappa$ will provide an explicit integration of our system.
Actually, we notice that both $|| \hat{S}_\pm ||$ are constants of the motion and $|| \hat{S}_+ ||^2 + || \hat{S}_- ||^2 = 2$, $\langle \hat{S}_+,\hat{S}_-\rangle = || \hat{S}_+ ||^2 - || \hat{S}_-||^2$.    The motion can be described by a precession around $B$ of the vector $\mathbf{S}_+$ of length $0 \leq || \mathbf{S}_+ || \leq 2$.  Once we have got the evolution of $\hat{S}_+$, then we can integrate the evolution equation for $\hat{S}_-$:
$$
\frac{d}{dt} \hat{S}_- = \mu [\hat{S}_-, \hat{B}] + \kappa [\hat{S}_-, \hat{S}_+] \, .
$$
It is remarkable that the quantity $\langle \hat{S}_1, \hat{S}_2 \rangle$ is a constant of the motion too (what can be checked after a simple computation), hence the system will not be state controllable, as the angle between the vectors determining the initial state will be preserved and not every configuration will be reachable.

A few simple solutions can be easily obtained.    For instance, it is obvious that if $\kappa = 0$ (this would happen, for instance, if $\mathbf{S}_1$ and $\mathbf{S}_2$ are parallel), then the equations describing the motion of the two systems, Eqs. (\ref{reducedhat})-(\ref{Pdothat}), decouple and each one behaves as the individual system described in Sect. \ref{sec:optimal} under the influence of the magnetic field $\mathbf{B}$, i.e., each spin precedes around the constant magnetic field $\mathbf{B}$.  

Notice, however, that the magnetic field is given by Eq. (\ref{optimalBhat}) and mixes both motions.    Thus, if we choose for instance $\mathbf{P}_1 = \mathbf{S}_2$ and  $\mathbf{P}_2 = -\mathbf{S}_1$, the magnetic field becomes $\mathbf{B} = 2\mu \mathbf{S_2} \wedge \mathbf{S}_1$ and the equations of motion become (using vector notation):
$$
\dot{\mathbf{S}}_1 = 2\mu^2 \, \mathbf{S}_1 \wedge (\mathbf{S}_2 \wedge \mathbf{S}_1 ) \, ; \qquad 
\dot{\mathbf{S}}_2 = 2\mu^2 \, \mathbf{S}_2 \wedge (\mathbf{S}_2 \wedge \mathbf{S}_1 ) \, .
$$
This system can be easily integrated using the identities for the triple vector product and noticing that $\cos \alpha = \mathbf{S}_1 \cdot \mathbf{S}_2$.   But, because $\mathbf{P}_1$ is perpendicular to $\mathbf{S}_1$, then $\cos\alpha = 0$. and the two perpendicular spins will rotate rigidily in the same plane.

Another family of solutions is obtained when the magnetic field vanishes or $\mu$ is negligible.  In such case the equations describing the motion become:
\begin{eqnarray*}
\frac{d}{dt}\hat{S}_1 &=& \kappa [\hat{S}_1, \hat{S}_2]\, , \qquad \frac{d}{dt}\hat{S}_2 =  \kappa  [\hat{S}_2,\hat{S}_1]  \\ 
\frac{d}{dt}\hat{P}_1 &=& \kappa [\hat{P}_1 - \hat{P}_2 , \hat{S}_2] \, , \qquad \frac{d}{dt}\hat{P}_2 = \kappa [\hat{P}_2 - \hat{P}_1,\hat{S}_1] \, .
\end{eqnarray*}
or, using the variables $\mathbf{S}_\pm$, we get from Eqs. (\ref{splusdot}), 
$$
\dot{\mathbf{S}}_+ = 0 \, , \qquad \dot{\mathbf{S}}_- = \kappa \mathbf{S}_- \wedge \mathbf{S}_+ \, ,
$$
or, in other words, $\mathbf{S}_-$ precedes around the constant vector $\mathbf{S}_+$.


\section{Conclusions}\label{sec:discussion}

The optimal control problem for two coupled spinning particles systems with given initial and final states, fixed time and control equations given in the Euler-Lagrange formalism, is analized and the differential equations determining their extremal solutions provided by PMP are exhibited. 

Such equations are obtained by using an adaptation of PMP for implicit control equations that involve a simple application of Dirac-Bergmann-Gotay constraints algorithm to reduce the implicit control differential equation of the system.  Thus the reduced control equations are obtained first and then PMP is applied.  

However it must be noticed that, in general, optimal control problems with implicit control equations doesn't necessarily satisfy this and, even for systems similar to the ones described in this paper, it could happen that the solutions to the optimal control problem are not found among the set of solution of the reduced control equations.  The space where the extremals are found being larger than the reduced space of the original implicit control equations.
This will happen, for instance, if the objective functional introduces some further constraints into the problem (for example like in the case of the time-optimal problem).  Under these more general circumstances another formulation of the problem is needed (like in the simple case of LQ systems discussed in \cite{De09}).
Results in this direction will be discussed in future publications.

It should also be pointed out that the resulting control problem for coupled spinning particles discussed in this paper excludes the possibility of entanglement, i.e., the reduced state space of the system is the Cartesian product of the corresponding individual reduced state spaces of the systems.   However, as pointed in the text, there are other couplings for which this is not true.  
Nevertherles, other possibilities could also be considered.  For instance, the group $SU(2) \times SU(2)$ is a diagonal subgroup of the group $SU(4)$ that could be considered as the configuration space of the composite system (in accordance with the quantum mechanical prescription and as commented already in the introduction \cite{Kh02}).   The exploration of such possibilities will be done elsewhere.

It was shown that the differential equations describing optimal extremals in the case of two coupled identical spinning particles in an uniform magnetic field with scalar coupling constitute a new, as far as it is known by the authors, completely integrable Hamiltonian system and a number of explicit solutions, both for individual and coupled systems, have been discussed.    The analysis performed here can be extended easily to chains of spinning particles a situation of interest in many applications as emphasized in the introduction.


\ack A.I. would like to thank P.J. Morrison  for pointing out that the equations describing the optimal control problem for a single classical spin, eqs. (\ref{SPmu}), are completely integrable.  R.S. would like to thank the support of project QUITEMAD, S-2009/ESP-1594. A.I. would also like to thank the organizers of the ``Second Iberoamerican Meeting on Geometry, Mechanics and Control'', where part of this work was completed.  

The authors would like to acknowledge the partial support by the Spanish MEC grant MTM2014-54692-P and QUITEMAD+, S2013/ICE-2801.  

We would also like to acknowledge the referees for their criticisms and valuable  suggestions.

\end{document}